\documentclass[fleqn]{mat01}
\usepackage{times,mathtimy,amssymb,latexsym}
\begin{document}

\setcounter{page}{319}
\firstpage{319}

\def\d{\hbox{d}}

\def\coth{\mbox{\rm coth}}

\newtheorem{theore}{Theorem}
\renewcommand\thetheore{\arabic{section}.\arabic{theore}}
\newtheorem{theor}[theore]{\bf Theorem}
\newtheorem{definit}[theore]{\rm DEFINITION}
\newtheorem{lem}[theore]{Lemma}
\newtheorem{rem}[theore]{Remark}
\newtheorem{propo}[theore]{\rm PROPOSITION}
\newtheorem{coro}[theore]{\rm COROLLARY}

\newcommand{\op}{\bigoplus}

\def\remar{\trivlist \item[\hskip \labelsep{\it Remarks.}]}

\title{Ideal amenability of Banach algebras on locally compact groups}

\markboth{M~Eshaghi Gordji and S~A~R~Hosseiniun}{Ideal amenability of Banach algebras}

\author{M~ESHAGHI GORDJI and S~A~R~HOSSEINIUN$^{*}$}

\address{Department of Mathematics, Faculty of Sciences, Semnan University, Semnan, Iran\\
\noindent $^{*}$Department of Mathematical Sciences, Shahid Beheshti University, Tehran, Iran\\
\noindent E-mail madjideg@walla.com; r-hosseini@cc.sbu.ac.ir}

\volume{115}

\mon{August}

\parts{3}

\pubyear{2005}

\Date{MS received 25 February 2005; revised 28 June 2005}

\begin{abstract}
In this paper we study the ideal amenability of Banach algebras. Let
$\cal A$ be a Banach algebra and let $I$ be a closed two-sided ideal
in $\cal A$, $\cal A$ is $I$-weakly amenable if $H^{1}({\cal
A},I^*)=\{0\}$. Further, $\cal A$ is ideally amenable if $\cal A$ is
$I$-weakly amenable for every closed two-sided ideal $I$ in $\cal
A$. We know that a continuous homomorphic image of an amenable
Banach algebra is again amenable. We show that for ideal amenability
the homomorphism property for suitable direct summands is true
similar to weak amenability and we apply this result for ideal
amenability of Banach algebras on locally compact groups.
\end{abstract}


\maketitle

\section{Introduction}

Let $A$ be a Banach algebra, then $A$ is amenable if
$H^1(A,X^{*})=\{0\}$ for each Banach $A$-module $X$; this definition was
introduced by Johnson in \cite{8}. The Banach algebra $A$ is weakly
amenable if $H^1(A,A^{*})=\{0\}$; this definition generalizes that
introduced by Bade, Curtis and Dales \cite{2}. Let $\cal A$ be a
Banach algebra and let $I$ be a closed two-sided ideal in $\cal A$,
$\cal A$ is $I$-weakly amenable if $H^{1}({\cal A},I^{*})=\{0\}$. ${\cal
A}$ is $n$-$I$-weakly amenable if $H^{1} ({\cal A}, I^{(n)}) = \{0\}$.

Further, $\cal A$ is ideally amenable if $\cal A$ is $I$-weakly amenable
for every closed two-sided ideal $I$ in $\cal A$ . This definition was
introduced by Eshaghi Gordji and Yazdanpanah in \cite{5}. In this paper
we study the ideal amenability of group algebras for a locally compact\break
group.

For example, it was shown in \cite{8} that the group algebra $L^1(G)$ is
amenable if and only if $G$ is an amenable group, and in \cite{9} that
$L^1(G)$ is weakly amenable for each locally compact group $G$. We know
that there exists a weakly amenable Banach algebra which is not ideally
amenable, also every $C^*$-algebra is ideally amenable \cite{5}. We can
not prove that for which locally compact group $G$, $L^1(G)$ is ideally
amenable. 

Let $A$ be a Banach algebra. The subspace $X$ of $A^{*}$ is left
introverted if $X^{*}\cdot X\subseteq X$ and $X^{*}$ by the following product
is a Banach algebra \cite{1}:
\begin{equation*}
\langle mn,f\rangle = \langle m,n\cdot f\rangle  \quad (m,n\in X^{*}, f\in X)
\end{equation*}
As in \cite{11} we have the following theorems.

\begin{theor}[\!]
Let $A$ be a Banach algebra and $Y\subseteq X$ two left introverted
subspaces of $A^{*}$. Then
\begin{equation*}
Y^\perp=\{m\in X^{*}\hbox{\rm :} \ \langle m,f \rangle = 0 \ \hbox{for every} \
f\in Y\}
\end{equation*}
is a $\hbox{weak}^{*}$ closed ideal of $X^{*}${\rm ,} also for each $n\in
X^{*}$ and $m\in Y^{\perp}$ we have $nm=0$.
\end{theor}

In this paper we consider $G$ to be a locally compact (topological)
group.

\begin{theor}[\!]
Let $X$ be a left introverted subspace of $L^\infty(G)$ such that
$C_{0}(G)\subseteq X\subseteq CB(G)$. Then there is an isometric
algebraic isomorphism $\tau\hbox{\rm :} \ M(G)\longrightarrow X^{*}$ and an isometric
direct sum decomposition $X^{*}=\tau(M(G))\oplus C_{0}(G)^\perp$ where
$C_{0}(G)^{\perp}$ is a $W^{*}$-closed ideal of $X^{*}$ and $X^{*}
C_{0}(G)^{\perp} = \{0\}$. 
\end{theor}

\begin{theor}[\!]
Let $N$ be a compact normal subgroup of a locally compact group $G${\rm ,} and let
$X_{N} = \{f\in C_{0}(G)\hbox{\rm :} \ l_{h} g = f \; for \; h\in N\}$.

Let $P\hbox{\rm :} \ C_{0}(G)\longrightarrow X_{N}$ be the following
linear map\hbox{\rm :}  
\begin{equation*}
\displaystyle (Pf)(x) = \int_{N} f(xh){\rm d} h,
\end{equation*} 
then $P^{*}\hbox{\rm :} \ (X_N)^{*}\longrightarrow M(G)$ is an isometric
isomorphism into $M(G)$ and there is the direct sum decomposition
$M(G)=P^{*}(X_{N})^{*}\oplus X_{N}^{\perp}$ where $P^{*}(X_{N})^{*}$ is
a closed subalgebra and $X_{N}^{\perp}$ is a weak$^*$ closed ideal of
$M(G)$ contained in right annihilators of $M(G)$.

Let $G$ be a locally compact group and $f\in CB(G)${\rm ,} then the right
orbit of $f$ is given by $R\circ(f)=\{r_xf\hbox{\rm :} \ x\in G\}$ and let
\begin{align*}
AP(G) &= \{f\in CB(G)\hbox{\rm :} \
R\circ (f) \; \hbox{is precompact in the norm topology}\},\\[.3pc]
W(G) &= \{f\in CB(G)\hbox{\rm :} \ R\circ (f) \; \hbox{is precompact in the
weak topology}\}.
\end{align*}
\end{theor}

\begin{theor}[\!]
There is a natural projection $Q\hbox{\rm :} \ W(G)\longrightarrow
AP(G)$ in which $Q^{*}\hbox{\rm :}$ $AP(G)^{*}\longrightarrow W(G)^{*}$
is an isometry and $W(G)^{*} = Q^{*}(AP(G)^{*})\oplus AP(G)^{\perp}$
where $Q^{*}(AP(G)^{*})$ is a unital closed subalgebra{\rm ;} and
$AP(G)^{\perp}$ is a weak$^*$ closed ideal of $W^{*}(G)${\rm ,} and
$AP(G)^{\perp}$ is contained in the right annihilators of $W(G)^{*}$.
\end{theor}

\begin{theor}[\!]
Let $G$ be amenable, and let $A\neq\{0\}$ be a $\hbox{weak}^{*}$ closed
self-adjoint translation invariant subalgebra of $L^{\infty}(G)$. Then
there exists a norm one projection $P\hbox{\rm :} \ L^{\infty}(G)\longrightarrow A$
which $L^{\infty}(G)^{*} = P^{*}(A^{*})\oplus A^{\perp}$ where
$A^{\perp}$ is a $\hbox{weak}^{*}$ closed ideal of $L^{\infty}(G)$ and
$P^{*}(A^{*})$ is a unital closed subalgebra of $L^{\infty}(G)^{*}$ and
$A^{\perp}$ is contained in the right annihilators of
$L^{\infty}(G)^{*}$. 
\end{theor}

As in \cite{4} for a left introverted subspace $X$ of $PM_{p}$ such that
$PF_{p}\subset X\subset UC_{p}$ we have as follows{\rm :}

\begin{theor}[\!]
There exists isometric algebra homomorphism $v\hbox{\rm :} \ W_{p}\rightarrow
X^{*}${\rm ,} such that $X^{*}=v(W_{p})\oplus PF_{p}^{\perp}.$\vspace{-.3pc}
\end{theor}

\section{Ideal amenability}

The concept of ideal amenability was introduced by Eshaghi Gordji and
Yazdanpanah \cite{5}. It is well-known that every ${\cal C}^*$-algebra
is ideally amenable. A commutative Banach algebra is ideally amenable if
and only if it is weakly amenable. Then ideal amenability is different 
from amenability. Let $\cal A$ be a non-unital non-weakly amenable
Banach algebra. Then there exists non-inner derivation $D\hbox{\rm :} \ {\cal A}
\rightarrow{\cal A}^*$. Let ${\cal A}^{\#}$ be the unitization of $\cal
A$, then we define $D_1\hbox{\rm :} \ {\cal A}^{\#} \rightarrow{\cal A}^*$ with $D_1
(a + \alpha) = D(a)$, $(a \in {\cal A}, \alpha \in \Bbb C)$. $D_1$ is a
derivation. Since $D$ is not inner, $D_1$ is not inner. ${\cal A}$ is a
closed two-sided ideal of ${\cal A}^{\#}$, therefore ${\cal A}^{\#}$ is
not ideally amenable. On the other hand, we know that there exists a 
non-weakly amenable Banach algebra $\cal A$ such that ${\cal A}^{\#}$ is
weakly amenable (see \cite{10}). Thus ${\cal A}^{\#}$ is an example of
weakly amenable Banach algebra which is not ideally amenable.

Suppose that $\cal A $ is a Banach algebra and let $X$ be $\cal
A$-bimodule. Let $ X {\cal A} = \{ f a\hbox{\rm :} \ f \in  X$, $a \in {\cal A} \}$
and ${\cal A} X = \{ a f\hbox{\rm :} \ a \in {\cal A}, f \in X \}.$ The space $X$
is said to be factors on the left (resp. right) if $X = X{\cal A} $
(resp. $ X={\cal A} X)$, and $X$ is neo-unital if $X$ is a factor on
both the left and the right. From Theorem~2.3 of \cite{7}, we have the
following.

\setcounter{theore}{0}

\begin{theor}[\!]
Let $\cal A$ be a Banach algebra with bounded approximate identity and
let $I$ be a codimension {\rm 1} closed two-sided ideal of $\cal A$. Then
$H^{1}({\cal A},X^{*})\cong H^{1}({\cal I},X^{*})$ for every neo-unital Banach
$\cal A$ bimodule $X$.
\end{theor}

\begin{coro} $\left.\right.$\vspace{.5pc}

\noindent Let $G$ be a discrete group{\rm ,} and let $I_0$ be a codimension {\rm 1}
closed two-sided ideal of $l^{1}(G)$. Then $l^{1}(G)$ is $I_0$-weakly
amenable.
\end{coro}

\begin{proof}
Since ${l^1(G)}$ is unital then ${I_0}$ is neo-unital
$l^{1}(G)$-bimodule. On the other hand, $I_0$ is weakly amenable
(Prop.~4.2 of \cite{7}). Then by the above theorem $l^{1}(G)$ is 
$I_0$-weakly amenable.

Let $X$, $Y$ and $Z$ be normed spaces and let $f\hbox{:}\ X\times
Y\longrightarrow Z$ be a continuous bilinear map, then
$f^{*}\hbox{:}\ Z^{*}\times X\longrightarrow Y^{*}$ (the transpose of $f$) is
defined by
\begin{equation*}
\langle f^{*}(z^{*}, x),y\rangle = \langle z^{*}, f(x, y)\rangle
\qquad (z^{*}\in Z^{*}, x\in X, y\in Y)
\end{equation*}
($f^*$ is a continuous bilinear map). Clearly, for each $x\in X$, the
mapping $z^*\longmapsto f^{*}(z^{*},x)\hbox{\rm :} \
Z^{*}\longrightarrow Y^*$ is $\hbox{weak}^{*}\hbox{--}\hbox{weak}^{*}$ continuous. Let
$f^{**}=(f^*)^*$ and $f^{***}=(f^{**})^*$, then $f^{***}\hbox{\rm :} \
X^{**}\times Y^{**}\longrightarrow Z^{**}$ is the unique extension of
$f$ such that $f^{***}(\cdot,y^{**})\hbox{\rm :} \ X^{**}\longrightarrow
Z^{**}$ is $\hbox{weak}^{*}\hbox{--}\hbox{weak}^{*}$ continuous. Also
$f^{***}(x,\cdot)\hbox{\rm :} \ Y^{**}\longrightarrow Z^{**}$ is
$\hbox{weak}^{*}\hbox{--}\hbox{weak}^{*}$ continuous for every $x\in X$. Let $
f^r\hbox{\rm :} \ Y\times X\longrightarrow Z$ be defined by $
f^r(y,x)=f(x,y)$, ~$(x\in X,~y\in Y)$. Then $f^r$ is a continuous linear
mapping from $Y\times X$ to $Z$. The map $f$ is Arens regular whenever
$f^{***}=f^{r***r}$. This is equivalent to the condition that the map
$f^{***}(x^{**},\cdot)\hbox{\rm :} \ Y^{**}\longrightarrow Z^{**}$ be
$\hbox{weak}^{*}\hbox{--}\hbox{weak}^{*}$ continuous for every $x^{**}\in X^{**}$. Let
$I$ be a closed two-sided ideal of $\cal A$ and $\pi_l\hbox{\rm :} \
{\cal A}\times I\longrightarrow I$ and $\pi_{r}\hbox{\rm :} \
I\times{\cal A}\longrightarrow I$ are the right and left module actions
of $\cal A$ on $I$ respectively. Then $I^{**}$ is a ${\cal
A}^{**}$-bimodule with module actions $\pi_l^{***}\hbox{\rm :} \ {\cal
A}^{**}\times I^{**}\longrightarrow I^{**}$ and $\pi^{***}_r\hbox{\rm :}
\ I^{**}\times{\cal A}^{**}\longrightarrow I^{**}$. $\pi_l$ is Arens
regular if and only if $a''$ is fixed in ${\cal A}^{**}$. The mapping
$i''\longmapsto \pi^{***}_l(i'',a'')\hbox{\rm :} \ I^{**}\longrightarrow
I^{**}$ is $\hbox{weak}^{*}\hbox{--}\hbox{weak}^{*}$ continuous and $\pi_r$ is Arens
regular if and only if $i''$ is fixed in $I^{**}$. The mapping
$a''\longmapsto \pi^{***}_r(a'',i'')\hbox{\rm :} \ {\cal
A}^{**}\longrightarrow I^{**}$ is $\hbox{weak}^{*}\hbox{--}\hbox{weak}^{*}$ continuous.
We now prove the following lemma.
\end{proof}

\begin{lem}
Let $\cal A$ be a Banach algebra with bounded approximate identity and
let $I$ be a closed two-sided ideal of $\cal A$. The the following
assertions hold.\vspace{-.4pc}
\begin{enumerate}
\renewcommand\labelenumi{\rm (\roman{enumi})}
\leftskip .2pc

\item If $\pi_r\hbox{\rm :} \ I\times{\cal A}\longrightarrow I$ be Arens regular{\rm ,} then
$I^*$ factors on the right.

\item If $\pi_l\hbox{\rm :} \ {\cal A}\times I\longrightarrow I$ be Arens regular{\rm ,}
then $I^*$ factors on the left.
\end{enumerate}
\end{lem}

\begin{proof} 
Let $(e_\alpha)$ be a bounded approximate identity for $\cal A$ with
cluster point $E$. Suppose that $\pi_r$ be Arens regular and let $i''\in
I^{**}$ be the cluster point of $(i_\beta)$ ($(i_\beta)$ is a net in
$I$). Then we have 
\begin{align*}
i'' &= \lim\limits_{\beta}\lim\limits_{\alpha} i_{\beta}\\[.3pc]
&= \lim\limits_{\beta}\lim\limits_{\alpha} i_{\beta} e_{\alpha}\\[.3pc]
&= \lim\limits_{\alpha}\lim\limits_{\beta} i_{\beta} e_{\alpha}\\[.3pc]
&= \lim\limits_{\alpha} i'' e_{\alpha} =i''E.
\end{align*}
Then for $i'\in I^*$ we have
\begin{align*}
\lim\limits_{\alpha}\langle e_{\alpha} i',i''\rangle &= \lim\limits_{\alpha}\langle i',
i''e_{\alpha}\rangle\\[.3pc]
&= \langle i', i''E \rangle\\[.3pc]
&= \langle i',i''\rangle.
\end{align*}
Thus $e_{\alpha} i' \rightarrow i'$ weakly. Since $e_{\alpha} i' \in
{\cal A}I^*$ for every $\alpha$, by the Cohen--Hewit factorization
theorem we know that ${\cal A}I^*$ is closed in $I^*$, then $i' \in
{\cal A}I^*$. Thus the proof of (i) is complete. For (ii), let $i''\in
I^{**}$ be the cluster point of $(i_\beta)$. For each $\beta$ we know
that $Ei_\beta=i_\beta$ since $\pi_l$ is Arens regular, then $e_{\alpha}
i'' \rightarrow Ei''=i''$ by $\hbox{weak}^*$ topology of $I^{**}$. Then for
every $i' \in I^*$,
\begin{align*}
\lim\limits_{\alpha}\langle i'e_{\alpha}, i''\rangle &= 
\lim\limits_{\alpha}\langle i', e_{\alpha} i''\rangle\\[.3pc]
&= \langle i', Ei'' \rangle\\[.3pc]
&= \langle i', i''\rangle.
\end{align*}
Therefore $i'e_{\alpha} \rightarrow i'$ weakly. Again by the Cohen--Hewit
factorization theorem we conclude that $I^*$ factors on the left.
\end{proof}

\begin{coro}$\left.\right.$\vspace{.5pc}

\noindent Let $\cal A$ be a Banach algebra with bounded approximate
identity and let $I$ be a codimension {\rm 1} closed two-sided ideal of
$\cal A$. If the module actions of $\cal A$ on $I$ are Arens regular{\rm ,}
then $\cal A$ is $2$-$I$-weakly amenable if and only if I is $2$-weakly
amenable.\vspace{-.3pc}
\end{coro}

\section{Ideal amenability and direct sum decompositions}

We prove an elementary lemma about derivations and homomorphic images
and apply this to ideal amenability for various algebras defined over
locally compact groups.

\setcounter{theore}{0}

\begin{lem}
Let $A$ be a right unital Banach algebra and $A=B\oplus I$ for closed
unital subalgebra $B$ and closed two-sided ideal $I$ in $A$ where $I$ is
contained in the set of right annihilators of $A$. Then we have the
following{\rm :} 
\begin{enumerate}
\renewcommand\labelenumi{\rm (\roman{enumi})}
\leftskip .2pc

\item for every closed two-sided ideal $J$ in $B${\rm ,} $J\oplus I$ is a closed
two-sided ideal in $A${\rm ,}

\item if $A$ is ideally amenable{\rm ,} then $B$ is ideally amenable.
\end{enumerate}
\end{lem}

\begin{proof} $\left.\right.$

\begin{enumerate}
\renewcommand\labelenumi{\rm (\roman{enumi})}
\leftskip .2pc

\item First we show that $J\oplus I$ is a two-sided ideal in $A$. Let $a
= b + i\in A$ and $x = j + i'\in J\oplus I$ since $AI=\{0\}$, then we
have
\begin{align*}
Ax &=(b+i)(j+i')=bj+ij\in J\oplus I,\\[.3pc]
xa &=(j+i')(b+i)=jb+i'b\in J\oplus I.
\end{align*}
Let $e$ be a right unit element in $A$, then there exists $b_1\in B$ and
$i_1\in I$ such that $e=b_1+i_1$. Then $e=e^{2}=e(b_{1}+i_{1})=eb_1$ and
$b_1$ is a right unit element for $A$. But $B$ is unital, and so
$b_1=1_B$. Let $\|1_B\|=M$ and $j+i\in J+I$, then we\break have 
\begin{equation*}
\|j+i\|\geq M^{-1}\| 1_B (j+i)\| = M^{-1}\|j\|.
\end{equation*}
Then the natural map $J\longrightarrow \frac{J+I}{I}$ is bicontinuous, therefore
$(J+I)/I$ is complete. $J$ is closed in ${A}/{I}$, and $J+I$ is
closed in $A$.

\item Let $J$ be a closed ideal in $B$ and $D\hbox{\rm :} \
B\longrightarrow J^*$ be a derivation and $\pi\hbox{\rm :} \ B\oplus
I\longrightarrow B$ be the natural embedding. $\pi$ is a homomorphism
and so is $\pi|_{J\oplus I}$. We show that $\bar{D}=(\pi|_{J\oplus
I})^*\circ D\circ\pi\hbox{\rm :} \ A\longrightarrow (J\oplus I)^*$ is a
derivation. For $x,y\in A$ and $z\in J\oplus I$, we have
\begin{align*}
\langle   (\pi|_{J\oplus I})^{*}\circ D\circ\pi (xy),z\rangle &=
\langle   D(\pi(x)\pi(y)), \pi(z)\rangle\\[.35pc]
&= \langle   D(\pi(x))\cdot  \pi(y)+\pi(x)\cdot  D(\pi(y)), \pi(z)\rangle,
\end{align*}
but
\begin{align*}
\langle   D(\pi(x))\cdot \pi(y), \pi(z)\rangle &=\langle   D(\pi(x)), \pi(y)
\pi(z)\rangle\\[.3pc]
&=\langle   D(\pi(x)), \pi(yz)\rangle\\[.35pc]
&=\langle   \bar{D}(x), yz\rangle\\[.35pc]
&=\langle   \bar{D}(x)\cdot y, z\rangle.
\end{align*}
Similarly,
\begin{equation*}
\langle \pi(x)\cdot D(\pi(y)), \pi(z)\rangle = \langle x\cdot
\bar{D}(y),z\rangle,
\end{equation*}
then $\bar{D}(xy)=x\cdot \bar{D}(y)+\bar{D}(x)\cdot y$ and $\bar{D}$ is
a derivation. $A$ is ideally amenable and by (i), $J\oplus I$ is a closed
two-sided ideal in $A$. Then there exists $\xi\in (J\oplus I)^*$ such that
$\bar{D}=\delta_\xi$, and let $\eta=\xi|_J$, then for $b\in B$ and $j\in
J$ we have
\begin{align*}
\langle Db,j\rangle &= \langle D(\pi(b)), \pi(j)\rangle\\[.35pc]
&=\langle  \bar{D} b, j\rangle\\[.35pc]
&=\langle b\xi-\xi\cdot b,j\rangle\\[.35pc]
&=\langle  \xi, j\cdot b\rangle  -\langle  \xi, b\cdot j\rangle\\[.35pc]
&=\langle b\cdot \eta-\eta\cdot b, j\rangle\\[.35pc]
&=\delta_\eta(b) (j).
\end{align*}
Therefore $D=\delta_\eta$ and $B$ is ideally amenable.\vspace{-1pc}
\end{enumerate}
\end{proof}\pagebreak

\begin{coro}$\left.\right.$\vspace{.5pc}

\noindent Let $\cal A$ be a dual Banach algebra with bounded approximate
identity and let ${\cal A}^{**}$ be the second dual of $\cal A$ with the
first Arens product. If ${\cal A}^{**}$ be ideally amenable then ${\cal
A}$ is ideally amenable.
\end{coro}

\begin{proof}
Let ${\cal A}_*$ be the pre-dual of $\cal A$. Then it is easy to see that
${\cal A}^{**}={\cal A}\oplus ({\cal A}_{*})^\perp$ (see, for example,
\cite{3} or \cite{6}). On the other hand ${\cal A}$ has bounded
approximate identity, then ${\cal A}^{**}$ has right identity. By lemma
the ideal amenability of ${\cal A}^{**}$ implies that ${\cal A}$ is also
ideally amenable.

We know that $M(G)$ is a dual Banach algebra (with pre-dual $C_0(G)$).
Thus by Coro\-llary~3.2, $M(G)$ is ideally amenable whenever ${M(G)}^{**}$ is
ideally amenable.
\end{proof}

\begin{coro}$\left.\right.$\vspace{.5pc}

\noindent Let $X$ be a left introverted subspace of $L^\infty(G)${\rm ,} such
that $C_0(G)\subseteq X\subseteq CB(G)${\rm ,} if $X^*$ be ideally amenable{\rm ,}
then $M(G)$ and $L^1(G)$ are ideally amenable.
\end{coro}

\begin{proof}
By applying Lemma~2.1 and Theorem~1.2, ideal amenability of $X^*$
implies that $M(G)$ is ideally amenable, and by \cite{5} we know that if
$M(G)$ is ideally amenable, then $L^1(G)$ is also ideally amenable.
\end{proof}

\begin{coro}$\left.\right.$\vspace{.5pc}

\noindent Let $N$ be a compact normal subgroup of a locally compact
group $G$, and $X_N$ and $P$ are as in Theorem~{\rm 1.3}. Then ideal
amenability of $M(G)$ implies that $P^*(X_N)^*$ is ideally amenable.
\end{coro}

\begin{coro}$\left.\right.$\vspace{.5pc}

\noindent In Theorem~{\rm 1.4,} if $W(G)^*$ be ideally amenable{\rm ,} then
$Q^{*}(AP(G)^{*})$ is ideally amenable.
\end{coro}

\begin{coro}$\left.\right.$\vspace{.5pc}

\noindent In Theorem~{\rm 1.5,} if $L^\infty(G)^*$ be ideally amenable{\rm ,} then
$P^*(A^*)$ is ideally amenable.
\end{coro}

Let $E$ be a right identity of $L^1(G)^{**}$ with $E\geq 0$ and
$\|E\|=1${\rm ,} then we have the isomorphism $EL^1(G)^{**}\simeq LUC (G)^*$
and we have the decomposition
\begin{align*}
L^1(G)^{**} &= EL^1(G)^{**}\oplus (I-E)L^1(G)^{**}\\[.3pc]
&= LUC (G)^*\oplus (I-E)L^1(G)^{**},
\end{align*}
where $(I-E)L^{1}(G)^{**}$ is continued in the right annihilators of
$L^{1}(G)^{**}$. So we have the decomposition $LUC (G)^{*}=M(G)\oplus
C_{0}(G)^\perp$ where $C_0(G)^\perp$ is a closed ideal in $LUC (G)^*$
and contained in the right annihilators of $LUC (G)^*$ (see, for
example, \cite{11}). Therefore the following proposition follows easily
from Theorem~2.1.

\begin{propo} $\left.\right.$\vspace{-.3pc}

\begin{enumerate}
\renewcommand\labelenumi{\rm (\roman{enumi})}
\leftskip .3pc

\item Suppose that $L^1(G)^{**}$ is ideally amenable{\rm ,} then $LUC (G)^*$ is
ideally amenable.

\item Suppose $LUC (G)^*$ is ideally amenable{\rm ,} then $M(G)$ is
ideally amenable.

\item Suppose $M(G)$ is ideally amenable{\rm ,} then $L^1(G)$ is ideally
amenable.
\end{enumerate}
\end{propo}

By Theorem~{\rm 1.6} we have as follows.\pagebreak

\begin{propo}$\left.\right.$\vspace{.5pc}

\noindent Let $X$ be a left introverted subspace of $PM_p${\rm ,} such that
$PM_p\subseteq X\subseteq UC_p${\rm ,} if $X^*$ be ideally amenable{\rm ,} then
$W_p$ is ideally amenable.
\end{propo}

\end{document}